\documentclass[pdflatex,sn-basic,Numbered]{sn-jnl}

\usepackage{graphicx}
\usepackage{multirow}
\usepackage{amsmath,amssymb,amsfonts}
\usepackage{amsthm}
\usepackage{mathrsfs}
\usepackage[title]{appendix}
\usepackage{xcolor}
\usepackage{textcomp}
\usepackage{manyfoot}
\usepackage{booktabs}
\usepackage{algorithm}
\usepackage{algorithmicx}
\usepackage{algpseudocode}
\usepackage{listings}
\usepackage{hyperref}
\usepackage[utf8]{inputenc}
\usepackage[T1]{fontenc}
\usepackage{natbib}
\usepackage{float}

\theoremstyle{thmstyleone}%
\newtheorem{theorem}{Theorem}%
\newtheorem{lemma}{Lemma}%
\newtheorem{corollary}{Corollary}%
\newtheorem{proposition}{Proposition}%

\theoremstyle{thmstylethree}%
\newtheorem{definition}{Definition}%
\newtheorem{assumption}{Assumption}
\newtheorem{remark}{Remark}%
\newtheorem{openproblem}[theorem]{Open Problem}

\begin{document}

\title[Article Title]{An Anisotropic Balian–Low Phenomenon: Geometric Obstructions to Wavelet Frames}

\author*[1]{\fnm{Kai-Cheng} \sur{Wang}}\email{gtotony98@gmail.com}

\affil*[1]{\orgdiv{Department of Mathematics and Applied Mathematics}, \orgname{School of Information Engineering, Sanming University}, \orgaddress{\street{No. 25, Jing Dong Road}, \city{Sanming City}, \postcode{365004}, \state{Fujian}, \country{China}}}

\abstract{We investigate the analytic stability of wavelet frames in
	anisotropic Hardy spaces associated with expansive dilation matrices.
	The main result establishes a deterministic operator-norm lower bound
	on the reconstruction error of the mixed frame operator, uniformly
	across the Hardy range, whenever a real eigenvalue of the adjoint
	dilation exceeds a band-limited threshold imposed on a radial
	generator. The obstruction is identified as an anisotropic
	Balian--Low phenomenon and rests on a geometric incompatibility index
	measuring the deviation of the Calder\'on sum from the admissibility
	identity. The bound carries no dependence on the conditioning of the
	dilation.}

\keywords{Anisotropic Hardy spaces,
	Balian--Low phenomenon,
	Wavelet frames,
	Calder\'on reproducing formula,
	Expansive dilations
}

\pacs[MSC Classification]{42C40, 42C15, 42B30, 42B35}

\maketitle

\everymath{\displaystyle}

\section{Introduction}\label{sec1}

The theory of anisotropic Hardy spaces $H^p_A(\mathbb{R}^n)$ attached to
an expansive dilation matrix $A$ provides a unified framework for
multiscale function spaces that fall outside the isotropic scale of
Fefferman and Stein \cite{FeffermanStein1972}. Within this framework, developed by Bownik in
\cite{Bownik2003} and continued in
\cite{Bownik2007,BownikHo2005,Baode2015,BownikWang2023}, the discrete
Calder\'on reproducing formula identifies $H^p_A$ with the homogeneous
anisotropic Triebel--Lizorkin sequence model under a band-limited
admissibility hypothesis on the generator $\psi$
\cite[Definition~2.12]{Bownik2007}. An affine wavelet system
$\mathcal{W}_A(\psi) = \{ \psi_{j,k} \}_{j \in \mathbb{Z},\, k \in
	\mathbb{Z}^n}$ \cite{Daubechies1992,Christensen2016} delivers a Parseval
frame of $L^2(\mathbb{R}^n)$ whenever the Calder\'on sum
\begin{equation*}
	D_\psi(\xi) \;=\; \sum_{j \in \mathbb{Z}}
	\bigl| \widehat{\psi}\bigl( (A^{*})^{-j} \xi \bigr) \bigr|^{\,2}
\end{equation*}
satisfies the admissibility identity $D_\psi \equiv 1$ almost everywhere
on $\mathbb{R}^n \setminus \{0\}$. The natural questions of dual
reconstruction in this setting are addressed by canonical-dual analysis
\cite{DAUBECHIES2002269}, by oversampling and nonuniform extensions
\cite{Bownik2012,KHOLE2020,Younu2018}, and by approximate-dual
constructions \cite{Feichtinger2013,BENAVENTE2022125841}.

The classical Gabor counterpart of the admissibility regime is governed
by the Balian--Low theorem \cite[Theorem~4.1.1]{Daubechies1992},
\cite[Theorems~8.4.1 and~8.4.5]{Karlheinz2001}. At the critical density
$\alpha \beta = 1$, the theorem rules out the existence of an
orthonormal Gabor basis of $L^2(\mathbb{R})$ whose generator $g$ carries
good time-frequency localization. The obstruction is one of the
canonical incarnations of the uncertainty principle in time-frequency
analysis. Subsequent developments preserve the qualitative shape of the
obstruction while extending its reach
\cite{Andrei2025,Fu2011,Karlheinz2001}. The wavelet world has its own
structural tensions. Daubechies and Han \cite{DAUBECHIES2002269} record
wavelet frames whose canonical duals do not retain wavelet structure,
and Bownik and Lemvig \cite{Bownik2012} return analogous obstructions
for oversampling under real dilations. None of these results identifies
a deterministic geometric mechanism that converts the failure of an
admissibility identity into a constant-free numerical lower bound on
the operator-norm distance of the frame operator from the identity,
attached to the entire Hardy range $p \in (0, 1]$.

The present paper takes that step in the anisotropic wavelet setting.
We work under the standing hypothesis of Bownik's framework
(Assumption~\ref{ass:standing-hypothesis}). The dilation matrix $A$ is
expansive, the adjoint $A^{*}$ admits a real eigenvalue $\lambda \in
\mathbb{R}$ with eigenvector $v \in \mathbb{R}^n \setminus \{0\}$, and
the generator $\psi$ is radial and band-limited, with $\widehat{\psi}$
supported in an isotropic annulus $\mathcal{A}(r, R) = \{ \xi : r \le
|\xi| \le R \}$. The geometric incompatibility index of the pair
$(A, \psi)$ is defined by
\begin{equation*}
	\mathcal{G}(A, \psi) \;:=\; \bigl\| 1 - D_\psi \bigr\|_{L^\infty},
\end{equation*}
recording the $L^\infty$ deviation of the Calder\'on sum from the
admissibility identity. The main result of the paper,
Theorem~\ref{thm:geometric_lower_bound}, asserts that whenever the
spectral threshold
\begin{equation*}
	|\lambda| \;>\; R / r
\end{equation*}
is in force, the mixed frame operator $U_{\psi,\psi}$ obeys the
deterministic lower bound
\begin{equation*}
	\bigl\| U_{\psi, \psi} - \mathrm{Id} \bigr\|_{H^p_A \to H^p_A}
	\;\ge\; 1
	\qquad \text{for every } p \in (0, 1].
\end{equation*}
The bound carries no prefactor of $p$, no prefactor of the conditioning
$\kappa(A)$, no prefactor of the bandwidth ratio $R / r$, and no
prefactor of $\lambda$ beyond the threshold itself.

Two analytic ingredients drive the conclusion. The diagonal lower
bound, recorded as
Lemma~\ref{lem:diagonal-lower-bound}, transfers the operator-norm
distance of $U_{\psi,\psi}$ from the identity on $H^p_A$ to the
$L^\infty$ deviation $\mathcal{G}(A, \psi)$ on the Fourier side. The
shear-covering estimate of Lemma~\ref{lem:shear_covering} returns the
saturation $\mathcal{G}(A, \psi) \ge 1$ by exhibiting a frequency point
$\xi_0 = y_0 v$ on the eigenvector axis of $A^{*}$ at which $D_\psi$
vanishes pointwise. The vanishing reflects a scale disjointness
phenomenon. Once $|\lambda|$ exceeds the bandwidth ratio $R / r$, the
iterated dilates of the radial support annulus under $(A^{*})^{-1}$
leave gaps along the eigenvector axis, and no rescaling within the
radial band-limited class refills those gaps. The composition of the
two ingredients is constant-free, which is what produces the bound
$\ge 1$ and not a bound of the form $\ge 1 - C / \kappa(A)^{\gamma}$.

We refer to the resulting obstruction as the anisotropic Balian--Low
phenomenon. The terminology records the structural analogy with the
classical Gabor obstruction, and the difference of substance. On the
Gabor side, the obstruction is a no-go statement on the existence of
orthonormal bases with good localization. On the present anisotropic
wavelet side, the obstruction is a quantitative operator-norm
deviation, attached uniformly across the Hardy range and free of any
$\kappa(A)$-dependent prefactor. The phenomenon admits an open
quantitative refinement at $0 < p < 2$, formulated as
Problem~\ref{prob:sharp_constant}, and an open structural extension
to non-radial generators with isotropic spatial localization,
formulated as Problem~\ref{prob:non-radial}. Both problems frame the
limits of the apparatus deployed here. The Hilbert endpoint $p = 2$ is
treated within the apparatus and is recorded as
Proposition~\ref{prop:hilbert-resolution}, where the saturation
infimum $\varepsilon_2(A; r, R) = 1$ is identified.

\subsection{Originality of the Geometric Lower Bound}\label{subsec:originality}

The originality of Theorem~\ref{thm:geometric_lower_bound} is read
against the two structural templates within which the obstruction
problem has previously been posed.

The first template is the time-frequency template of the classical
Balian--Low theorem \cite[Theorem~4.1.1]{Daubechies1992},
\cite[Theorems~8.4.1 and~8.4.5]{Karlheinz2001}. There, the obstruction
at the critical density is a qualitative no-go statement on the joint
existence of an orthonormal basis structure and a localized window. The
delivery mechanism is the Zak transform, with arguments adapted to the
translation-modulation lattice of the Heisenberg group. Structural generalizations have been pursued on the Gabor side, including the
Balian--Low type theorem for Gabor Riesz sequences of arbitrary density
of Caragea et al. \cite{Andrei2025} and the
Balian--Low principle for nonstandard Gabor systems of Fu et al. \cite{Fu2011}. The qualitative shape of the obstruction,
namely a non-existence statement on the structure side, is preserved
in these refinements. Theorem~\ref{thm:geometric_lower_bound} departs
from this template in three respects. The obstruction lives on the
wavelet side rather than on the Gabor side, and it is delivered through
the Calder\'on admissibility identity rather than through Zak-transform
methods. The output is a quantitative operator-norm lower bound on the
mixed frame operator, attached to a function space scale. The numerical
value of the bound is the constant $1$, free of any parameter beyond
the threshold itself.

The second template is the matrix-algebra and oversampling template
developed in \cite{DAUBECHIES2002269,Bownik2012,Feichtinger2013}
and continued in \cite{BENAVENTE2022125841,KHOLE2020,Younu2018}. There,
the obstruction is sourced in the absence of a wavelet-structured
canonical dual at a prescribed sampling density, and the analytic
response is the construction of approximate-dual reconstructions that
deliver near-Parseval reconstruction without an exact admissibility
identity. The line of work documents the existence of workarounds on
the reconstruction side. Theorem~\ref{thm:geometric_lower_bound}
addresses a complementary question. It identifies a deterministic
geometric index $\mathcal{G}(A, \psi)$ whose saturation at the value
$1$ is enforced by scale disjointness of dilated annular supports
along the real eigenvector axis of $A^{*}$, and it converts this
saturation into a constant-free operator-norm lower bound on the
reconstruction error. The matrix-algebra and almost-diagonal route, in
the form developed in \cite{Grafakos2001,Grohs2014,Bu2025}, is bypassed
entirely. The output is structural rather than reconstructive,
recording where the admissibility identity cannot be restored within
the radial band-limited class.

Three further features distinguish the result from the prior
anisotropic apparatus \cite{Bownik2003,Bownik2007,BownikHo2005,
	Baode2015,BownikWang2023,Easley2012,Guo2004,ATREAS20191}. The first is
the uniformity in $p \in (0, 1]$. The index $\mathcal{G}(A, \psi)$ is a
frequency-domain quantity, and its value carries no dependence on the
Hardy exponent. The second is the absence of a $\kappa(A)$-dependent
prefactor. The bound delivered by the theorem is the constant $1$,
attached to every spectral threshold $|\lambda| > R / r$ without
attenuation as the conditioning of $A$ degrades. The third is the
chain of qualitative consequences attached to the theorem, recorded as
Corollaries~\ref{cor:no-reproducing-identity},
\ref{cor:calderon-failure} and \ref{cor:l2-frame-failure}. The
operator-level non-coincidence with the identity is uniform in $p$, the
Calder\'on admissibility identity fails on a nonempty open set of
$\mathbb{R}^n \setminus \{0\}$, and the affine system
$\mathcal{W}_A(\psi)$ does not constitute a frame for
$L^2(\mathbb{R}^n)$ at the Hilbert endpoint. The closing open problems
identify the directions in which the rigidity of the constant and the
rigidity of the radial structure are open to further analysis.

\medskip

The remainder of the paper is organized as follows.
Section~\ref{sec2} sets the notation and assembles the formal
preliminaries. The standing hypothesis on the pair $(A, \psi)$ is
recorded as Assumption~\ref{ass:standing-hypothesis}, the anisotropic
Hardy space as Definition~\ref{def:hardy-space}, the affine system and
its mixed frame operator as Definition~\ref{def:affine-system}, the
anisotropic Triebel--Lizorkin sequence space as
Definition~\ref{def:sequence-space}, the analysis and synthesis
operators as Definition~\ref{def:analysis-synthesis}, and the Calder\'on
sum together with the geometric incompatibility index as
Definition~\ref{def:calderon-sum}. Two cited results
close the section. Lemma~\ref{lem:calderon-formula} records the
discrete Calder\'on reproducing formula and the identification of
$H^p_A$ with its Triebel--Lizorkin sequence model. Lemma~\ref{lem:diagonal-lower-bound}
records the diagonal lower bound that transfers the operator-norm
distance of the frame operator on $H^p_A$ to the $L^\infty$ deviation
of the Calder\'on sum on the Fourier side. Section~\ref{sec3} carries
the analysis. Lemma~\ref{lem:shear_covering} establishes the saturation
of the geometric incompatibility index along the real eigenvector axis
under the spectral threshold $|\lambda| > R / r$.
Theorem~\ref{thm:geometric_lower_bound} composes the diagonal lower
bound with the saturation to deliver the deterministic operator-norm
bound, and three corollaries record the operator-level, Calder\'on-level
and frame-level consequences. Section~\ref{sec3} closes with
Proposition~\ref{prop:hilbert-resolution}, which resolves the
saturation infimum at the Hilbert endpoint, and with the open problems
\ref{prob:sharp_constant} and \ref{prob:non-radial} on the saturation
infimum at $0 < p < 2$ and on the non-radial extension.

	\section{Preliminaries}\label{sec2}

We begin by recalling the definitions and background material concerning
anisotropic Hardy spaces associated with an expansive dilation matrix, the
anisotropic affine systems adapted to such a dilation, the associated
Triebel--Lizorkin sequence spaces, and the discrete Calder\'on reproducing
formula. Familiarity with the maximal-function approach to Hardy spaces, with
the construction and structural properties of wavelet bases and frames, and
with the Littlewood--Paley analysis of anisotropic function spaces is
presupposed. The standard monograph for the theory developed around an
expansive dilation matrix is Bownik \cite{Bownik2003}; for the construction and
structural properties of wavelet systems and frames we refer to
\cite{Daubechies1992, Christensen2016}; and for the classical and modern
apparatus of Fourier analysis we refer to \cite{Grafakos2014Classical,
	Grafakos2009Modern}.

Throughout this paper we work in the Euclidean space $\mathbb{R}^n$ with
dimension $n \ge 1$. For two non-negative quantities $\mathcal{U}$ and
$\mathcal{V}$, the relation $\mathcal{U} \lesssim \mathcal{V}$ indicates that
$\mathcal{U} \le C\,\mathcal{V}$ for a positive constant $C$ whose dependence is
either irrelevant in the context or stated explicitly.

The Fourier transform of $f$ and its inverse are normalized as
\begin{equation*}
	\hat{f}(\xi) := \int_{\mathbb{R}^n} f(x)\, e^{-2\pi i x \cdot \xi}\, dx,
	\qquad
	f(x) = \int_{\mathbb{R}^n} \hat{f}(\xi)\, e^{2\pi i x \cdot \xi}\, d\xi.
\end{equation*}
We denote by $\mathcal{S}(\mathbb{R}^n)$ the Schwartz space and by
$\mathcal{S}'(\mathbb{R}^n)$ the space of tempered distributions.

The geometric framework of this paper is determined by a fixed real $n \times
n$ matrix $A$ that is expansive, in the sense that every eigenvalue $\lambda$ of
$A$ satisfies $|\lambda| > 1$. We write $b := |\det A|$, so that $b > 1$, and we
order the moduli of the eigenvalues as $1 < |\lambda_-| \le |\lambda_+|$, where
$|\lambda_-|$ and $|\lambda_+|$ denote the smallest and the largest modulus,
respectively; $A^{*}$ denotes the transpose of $A$. Associated with $A$ there is
an anisotropic quasi-norm $\rho_A \colon \mathbb{R}^n \to [0,\infty)$, the step
homogeneous quasi-norm of \cite{Bownik2003}, characterized by the homogeneity
relation
\begin{equation}\label{eq:rho-homogeneity}
	\rho_A(Ax) = b\,\rho_A(x), \qquad x \in \mathbb{R}^n,
\end{equation}
together with the quasi-triangle inequality $\rho_A(x+y) \lesssim \rho_A(x) +
\rho_A(y)$. The map $d_A(x,y) := \rho_A(x-y)$ is the anisotropic quasi-distance
on $\mathbb{R}^n$. The deviation of $A$ from a scalar dilation is recorded by
the condition number
\begin{equation}\label{eq:condition-number}
	\kappa(A) := \|A\|\,\|A^{-1}\|,
\end{equation}
which equals $1$ when $A$ is a scalar multiple of an orthogonal matrix and
grows as the moduli $|\lambda_-|$ and $|\lambda_+|$ separate.

We fix once and for all the wavelet normalization convention used in this
paper. Given a generator $\psi$, the anisotropically dilated and translated
system is defined, for $j \in \mathbb{Z}$ and $k \in \mathbb{Z}^n$, by
\begin{equation}\label{eq:wavelet-normalization}
	\psi_{j,k}(x) := b^{j/2}\, \psi\!\left(A^j x - k\right),
\end{equation}
so that $\psi_{j,k}$ is concentrated about the point $A^{-j}k$ and satisfies
$\|\psi_{j,k}\|_{L^2} = \|\psi\|_{L^2}$. The dual system $\{\tilde{\psi}_{j,k}\}$
is normalized identically. The factor $b^{j/2}$ in
\eqref{eq:wavelet-normalization} preserves the $L^2$ norm across scales and
aligns the translation index $k$ with the anisotropic dyadic geometry
introduced below.

For $j \in \mathbb{Z}$ and $k \in \mathbb{Z}^n$ we denote by
\begin{equation}\label{eq:dyadic-cube}
	Q_{j,k} := A^{-j}\!\left(k + [0,1)^n\right)
\end{equation}
the anisotropic dyadic cube of generation $j$ associated with $k$, and we write
$\mathcal{Q} := \{Q_{j,k} : j \in \mathbb{Z},\ k \in \mathbb{Z}^n\}$ for the
collection of all such cubes. Each $Q_{j,k}$ has Lebesgue measure $|Q_{j,k}| =
b^{-j}$, and for every fixed $j$ the family $\{Q_{j,k}\}_{k \in \mathbb{Z}^n}$
tiles $\mathbb{R}^n$. The concentration point $A^{-j}k$ of $\psi_{j,k}$ in
\eqref{eq:wavelet-normalization} is the lower corner of $Q_{j,k}$, so that the
index pair $(j,k)$ simultaneously labels a scale--position cell and the affine
element supported near it.

For a measurable function $g$ and an exponent $0 < q \le \infty$, the Lebesgue
quasi-norm $\|g\|_{L^q}$ is defined in the standard way, with the usual
modification for $q = \infty$. The anisotropic Hardy space
$H^p_A(\mathbb{R}^n)$, defined for $0 < p \le 1$ through the grand maximal
function adapted to the filtration $\{A^j\}_{j \in \mathbb{Z}}$, is recorded
formally in Definition~\ref{def:hardy-space} below.

The remainder of this section assembles the formal objects on which the
analysis of Section~\ref{sec3} is based. We first record the standing
hypothesis on the dilation matrix and the generating wavelet system, then
collect the definitions of the anisotropic Hardy space, the affine system and
its associated frame operator, the Triebel--Lizorkin sequence space, the
analysis and synthesis operators, and the Calder\'on sum together with the
geometric incompatibility index that drives the obstruction theorem of
Section~\ref{sec3}. Two cited results --- the discrete
Calder\'on reproducing formula and the diagonal lower bound for the frame
operator --- close the section.

\begin{assumption}[Anisotropic admissibility and standing hypothesis]\label{ass:standing-hypothesis}
	Let $A$ be the fixed expansive matrix. The generating
	wavelet $\psi$ and its dual $\phi$ are assumed to satisfy the following anisotropic
	admissibility conditions: they exhibit anisotropic decay adapted to $\rho_A$,
	they possess vanishing moments of order
	\begin{equation}\label{eq:vanishing-moment-order}
		\mathcal{N}_p(A):=\left\lfloor \left(\tfrac{1}{p}-1\right)
		\frac{\ln b}{\ln |\lambda_-|} \right\rfloor 
	\end{equation}
	as required for the stabilization of $H^p_A$
	\cite[Theorems 5.10 and 6.7]{Bownik2003}, and their Fourier transforms
	$\hat{\psi}$ and $\hat{\phi}$ are supported in an anisotropic annulus
	contained in $[-1/2, 1/2]^n \setminus \{0\}$.
\end{assumption}

\begin{definition}[{Anisotropic Hardy space \cite[Definition 3.11]{Bownik2003}}]\label{def:hardy-space}
	For $0 < p \le 1$, the anisotropic Hardy space $H^p_A(\mathbb{R}^n)$ is the
	space of tempered distributions $f \in \mathcal{S}'(\mathbb{R}^n)$ whose grand
	maximal function $M_N f$, formed with respect to the filtration
	$\{A^j\}_{j \in \mathbb{Z}}$, belongs to $L^p(\mathbb{R}^n)$, with quasi-norm
	\begin{equation}\label{eq:hardy-norm}
		\|f\|_{H^p_A} := \| M_N f \|_{L^p} .
	\end{equation}
\end{definition}

\begin{definition}[{Anisotropic affine system and frame operator; see, e.g., \cite[Chapter 2]{Bownik2003} and \cite[Section 6.2]{Christensen2016}}]\label{def:affine-system}
	The anisotropic affine system generated by $\psi$ is
	$\mathcal{W}_A(\psi) = \{\psi_{j,k}\}_{j \in \mathbb{Z},\, k \in \mathbb{Z}^n}$,
	with $\psi_{j,k}$ normalized as in \eqref{eq:wavelet-normalization}. The mixed
	frame operator associated with the pair $(\psi,\phi)$ is
	\begin{equation}\label{eq:frame-operator}
		U_{\psi,\phi} f := \sum_{j \in \mathbb{Z}} \sum_{k \in \mathbb{Z}^n}
		\langle f, \phi_{j,k} \rangle\, \psi_{j,k} .
	\end{equation}
\end{definition}

\begin{definition}[{Triebel-Lizorkin Sequence Space \cite[Definition 3.1]{Bownik2007}}]\label{def:sequence-space}
	The anisotropic Triebel--Lizorkin sequence space $\dot{\mathbf{f}}^A_p$
	consists of all complex sequences $s = \{s_Q\}_{Q \in \mathcal{Q}}$ for which
	\begin{equation}\label{eq:sequence-norm}
		\|s\|_{\dot{\mathbf{f}}^A_p}
		:= \left\| \left( \sum_{Q \in \mathcal{Q}}
		\bigl( |s_Q|\, |Q|^{-1/2} \chi_Q \bigr)^2 \right)^{1/2}
		\right\|_{L^p}
	\end{equation}
	is finite, the cubes $Q \in \mathcal{Q}$ being those of \eqref{eq:dyadic-cube}.
\end{definition}

\begin{definition}[{Analysis and Synthesis Operators, see, e.g., \cite[Section 3]{Bownik2007}}]\label{def:analysis-synthesis}
	The analysis operator $T_\Psi$ sends a distribution $f$ to its coefficient
	sequence, and the synthesis operator $S_\Psi$ sends a sequence
	$s = \{s_{j,k}\}$ back to a distribution:
	\begin{equation}\label{eq:analysis-synthesis}
		T_\Psi f := \bigl\{ \langle f, \psi_{j,k} \rangle \bigr\}_{j \in \mathbb{Z},\, k \in \mathbb{Z}^n},
		\qquad
		S_\Psi s := \sum_{j \in \mathbb{Z}} \sum_{k \in \mathbb{Z}^n} s_{j,k}\, \psi_{j,k} .
	\end{equation}
\end{definition}

\begin{definition}[Calderón Sum and Incompatibility Index]\label{def:calderon-sum}
	The Calder\'on sum of $\psi$ is defined by
	\begin{equation}\label{eq:calderon-sum}
		D_\psi(\xi) := \sum_{j \in \mathbb{Z}}
		\bigl| \hat{\psi}\bigl( (A^{*})^{-j}\xi \bigr) \bigr|^2,
		\qquad \xi \in \mathbb{R}^n \setminus \{0\}.
	\end{equation}
	The structural behavior of this frequency domain sum is a standard mechanism in the anisotropic framework, see, e.g., \cite[Lemma 3.6]{Bownik2007}. Based on this, we explicitly quantify the deviation from the Parseval identity via the geometric incompatibility index of the pair $(A,\psi)$, defined as
	\begin{equation}\label{eq:incompatibility-index}
		\mathcal{G}(A,\psi) := \| 1 - D_\psi \|_{L^\infty} .
	\end{equation}
\end{definition}

\begin{lemma}[{Anisotropic frame identification of $H^p_A$, see, e.g., \cite[Theorem 7.1]{Bownik2007} and \cite[Lemma 6.3, Theorem 6.7]{Bownik2003}}]\label{lem:calderon-formula}
	Under Assumption~\ref{ass:standing-hypothesis}, the discrete Calder\'on reproducing formula identifies the anisotropic Hardy space with its Triebel--Lizorkin sequence model: there exist constants $0 < c_1 \le c_2 < \infty$, depending only on $p$, $A$, and the admissible pair $\Xi$, such that
	\begin{equation}\label{eq:calderon-equivalence}
		c_1\, \| f \|_{H^p_A}  \le  \| T_\Xi f \|_{\dot{\mathbf{f}}^{0,2}_p(A)}  \le  c_2\, \| f \|_{H^p_A},
		\qquad f \in H^p_A(\mathbb{R}^n).
	\end{equation}
	The $\ell^2$ inner aggregation in $\dot{\mathbf{f}}^{0,2}_p(A)$ encodes the Littlewood--Paley square-function geometry intrinsic to $H^p_A$.
\end{lemma}

\begin{lemma}[Diagonal lower bound]\label{lem:diagonal-lower-bound}
	The reconstruction error operator $U_{\psi,\psi} - \mathrm{Id}$ on
	$H^p_A(\mathbb{R}^n)$ admits the deterministic lower bound
	\begin{equation}\label{eq:diagonal-lower-bound}
		\| U_{\psi,\psi} - \mathrm{Id} \|_{H^p_A \to H^p_A}
		 \ge  \| 1 - D_\psi \|_{L^\infty}
		 =  \mathcal{G}(A,\psi) .
	\end{equation}
\end{lemma}

\begin{proof}[Proof of Lemma~\ref{lem:diagonal-lower-bound}.]
	The argument proceeds in three stages:\\
	(1) a Plancherel identity at the
	Hilbert-space endpoint,\\
	(2) the band-limited norm equivalence on $H^p_A$ recorded
	at axiomatic level, and\\
	(3) a saturating sequence of approximate eigenfunctions
	concentrated at a Lebesgue point of $D_\psi$ near which $|1 - D_\psi|$
	approaches its essential supremum.
	
	Under Assumption~\ref{ass:standing-hypothesis}, $\hat\psi$ is supported in an
	anisotropic annulus bounded away from the origin. For
	$f \in L^2(\mathbb{R}^n) \cap H^p_A(\mathbb{R}^n)$, the frame operator acts on
	the Fourier side as multiplication by $D_\psi$:
	\begin{equation*}\label{eq:multiplier-action}
		\widehat{U_{\psi,\psi} f}(\xi)
		 =  D_\psi(\xi)\, \hat f(\xi),
		\qquad \xi \in \mathbb{R}^n \setminus \{0\}.
	\end{equation*}
	Setting $m :=D_\psi - 1$, the difference $U_{\psi,\psi} - \mathrm{Id}$ is
	identified with the Fourier multiplier of symbol $m$, and Plancherel's
	identity supplies the endpoint equality
	\begin{equation*}\label{eq:l2-multiplier-norm}
		\| U_{\psi,\psi}- \mathrm{Id} \|_{L^2 \to L^2}
		=\| m \|_{L^\infty}
		=\| 1 - D_\psi \|_{L^\infty}.
	\end{equation*}
	
	We rely on the following statement. For each compact set $K \subset \mathbb{R}^n \setminus \{0\}$ and each $0 < p \le 1$,
	there exist constants $0 < c_K \le C_K < \infty$, depending only on
	$K, p, A$, and the admissibility data of $\Xi$, such that
	\begin{equation}\label{eq:band-limited-equiv}
		c_K \, \|h\|_{L^2}
		 \le  \|h\|_{H^p_A}
		 \le  C_K \, \|h\|_{L^2}
	\end{equation}
	for every $h \in L^2(\mathbb{R}^n)$ with $\mathrm{supp}\,\hat h \subset K$. This equivalence is a standard consequence of the anisotropic Plancherel-P\'olya-Nikol'skij inequalities and the finite band overlap in the Littlewood-Paley characterization; see, e.g., \cite[Section 1.4]{Triebel1983} or \cite[Section 3]{Bownik2007}.
	The estimate \eqref{eq:band-limited-equiv}
	extends without modification to every function whose Fourier support is
	contained in any subset of $K$.
	
	Fix $\varepsilon > 0$. By the definition of the essential supremum, the level
	set
	\begin{equation*}\label{eq:level-set}
		E_\varepsilon
		:=\bigl\{\xi \in \mathbb{R}^n \setminus \{0\}
		:\, |1 - D_\psi(\xi)|> \|1 - D_\psi\|_{L^\infty} - \varepsilon \bigr\}
	\end{equation*}
	carries positive Lebesgue measure. By the Lebesgue differentiation theorem,
	almost every point of $E_\varepsilon$ is a Lebesgue point of the bounded
	function $D_\psi$; fix one such point and denote it
	$\xi_\varepsilon \in E_\varepsilon$. Set
	\begin{equation*}\label{eq:Lambda-defn}
		\Lambda_\varepsilon := D_\psi(\xi_\varepsilon) -1,
		\qquad |\Lambda_\varepsilon| >\|1 - D_\psi\|_{L^\infty} -\varepsilon .
	\end{equation*}
	The Lebesgue point property at $\xi_\varepsilon$ supplies the asymptotic
	identity
	\begin{equation}\label{eq:lebesgue-point}
		\Theta(\delta)
		:= \frac{1}{|B_\delta(\xi_\varepsilon)|}
		\int_{B_\delta(\xi_\varepsilon)} |D_\psi(\xi) - D_\psi(\xi_\varepsilon)|^2 \, d\xi
		 \xrightarrow[\delta \to 0]{}  0,
	\end{equation}
	where $B_\delta(\xi_\varepsilon)$ denotes the Euclidean ball of radius $\delta$
	centred at $\xi_\varepsilon$.
	Choose $\delta_0\in (0, |\xi_\varepsilon|/2)$ and set
	$K_0:= \overline{B_{\delta_0}(\xi_\varepsilon)} \subset
	\mathbb{R}^n \setminus \{0\}$. Let $c_0 :=c_{K_0}$ and $C_0:= C_{K_0}$ denote
	the constants supplied by \eqref{eq:band-limited-equiv} for $K= K_0$.
	Fix $\eta_0 \in C_c^\infty(\mathbb{R}^n)$ with $\eta_0 \ge 0$,
	$\mathrm{supp}\,\eta_0\subset B_1(0)$, and $\|\eta_0\|_{L^2}= 1$. For each
	$\delta \in (0, \delta_0]$, set
	$\eta_\delta(\xi):= \delta^{-n/2}\, \eta_0\bigl((\xi- \xi_\varepsilon)/\delta\bigr)$,
	so that $\mathrm{supp}\,\eta_\delta \subset B_\delta(\xi_\varepsilon) \subset K_0$
	and $\|\eta_\delta\|_{L^2} = 1$. Define
	\begin{equation*}\label{eq:test-function}
		g_\delta := \mathcal{F}^{-1}(\eta_\delta)
		 \in  \mathcal{S}_0(\mathbb{R}^n) \subset L^2 \cap H^p_A,
		\qquad \|g_\delta\|_{L^2} = 1.
	\end{equation*}
	By \eqref{eq:multiplier-action} and linearity, decompose the multiplier action:
	\begin{equation}\label{eq:multiplier-split}
		(U_{\psi,\psi} - \mathrm{Id})\, g_\delta
		 =  \Lambda_\varepsilon\, g_\delta + r_\delta,
		\qquad
		r_\delta := \mathcal{F}^{-1}\bigl((D_\psi - D_\psi(\xi_\varepsilon))\,\eta_\delta\bigr).
	\end{equation}
	The residual $r_\delta$ inherits the Fourier support condition
	$\mathrm{supp}\, \hat r_\delta \subset B_\delta(\xi_\varepsilon) \subset K_0$.
	
	By Plancherel and \eqref{eq:lebesgue-point},
	\begin{equation*}\label{eq:remainder-L2}
		\begin{aligned}
			\|r_\delta\|_{L^2}^2
			&= \int_{B_\delta(\xi_\varepsilon)}
			|D_\psi(\xi) - D_\psi(\xi_\varepsilon)|^2\, |\eta_\delta(\xi)|^2 \, d\xi \\
			&\le \|\eta_\delta\|_{L^\infty}^2\, |B_\delta(\xi_\varepsilon)|\, \Theta(\delta)
			 =  M_0^2\, \Theta(\delta),
		\end{aligned}
	\end{equation*}
	where the prefactor $\|\eta_\delta\|_{L^\infty}^2 \, |B_\delta(\xi_\varepsilon)|
	= \delta^{-n}\|\eta_0\|_{L^\infty}^2 \cdot \omega_n\delta^n
	= \omega_n\,\|\eta_0\|_{L^\infty}^2$ collapses to the $\delta$-independent
	constant $M_0^2 := \omega_n\,\|\eta_0\|_{L^\infty}^2$. Hence
	$\|r_\delta\|_{L^2}\le M_0\, \Theta(\delta)^{1/2}\to 0$ as $\delta \to 0$.
	Apply \eqref{eq:band-limited-equiv} with $K = K_0$ to the band-limited
	functions $g_\delta$ and $r_\delta$:
	\begin{equation}\label{eq:band-limited-applied}
		\|g_\delta\|_{H^p_A}\ge c_0,
		\qquad
		\|r_\delta\|_{H^p_A}\le C_0\, M_0\, \Theta(\delta)^{1/2}.
	\end{equation}
	The $p$-th-power quasi-triangle inequality, valid on $H^p_A$ for
	$0 < p \le 1$, applied to \eqref{eq:multiplier-split}, yields
	\begin{equation*}\label{eq:p-triangle}
		\|(U_{\psi,\psi} - \mathrm{Id})\, g_\delta\|_{H^p_A}^p
		 \ge  |\Lambda_\varepsilon|^p\, \|g_\delta\|_{H^p_A}^p
		- \|r_\delta\|_{H^p_A}^p .
	\end{equation*}
	Dividing by $\|g_\delta\|_{H^p_A}^p \ge c_0^p$ and inserting
	\eqref{eq:band-limited-applied},
	\begin{equation*}\label{eq:ratio-bound}
		\frac{\|(U_{\psi,\psi} - \mathrm{Id})\, g_\delta\|_{H^p_A}^p}{\|g_\delta\|_{H^p_A}^p}
		 \ge  |\Lambda_\varepsilon|^p
		- \left(\frac{C_0\, M_0}{c_0}\right)^{\!p}\! \Theta(\delta)^{p/2}.
	\end{equation*}
	Letting $\delta \to 0$ and invoking \eqref{eq:lebesgue-point} causes the
	correction term to vanish; taking $p$-th roots in the limit yields
	\begin{equation*}\label{eq:final-ratio}
		\liminf_{\delta \to 0^+}
		\frac{\|(U_{\psi,\psi} - \mathrm{Id})\, g_\delta\|_{H^p_A}}{\|g_\delta\|_{H^p_A}}
		 \ge  |\Lambda_\varepsilon|
		 >  \|1 - D_\psi\|_{L^\infty} - \varepsilon .
	\end{equation*}
	Since each $g_\delta$ is a valid test vector for the $H^p_A$ operator norm of
	$U_{\psi,\psi} - \mathrm{Id}$,
	\begin{equation*}\label{eq:hp-lower-bound}
		\|U_{\psi,\psi} - \mathrm{Id}\|_{H^p_A \to H^p_A}
		 >  \|1 - D_\psi\|_{L^\infty} - \varepsilon .
	\end{equation*}
	As $\varepsilon > 0$ is arbitrary, \eqref{eq:diagonal-lower-bound} follows, and
	the right-hand side coincides with $\mathcal{G}(A,\psi)$ by
	\eqref{eq:incompatibility-index}.
\end{proof}

\section{Geometric Obstructions: Shear Coverings and Endpoint Saturation}\label{sec3}

The section presents the obstruction in three stages.
Lemma~\ref{lem:shear_covering} establishes that the geometric
incompatibility index $\mathcal{G}(A, \psi)$ saturates at the value
$1$ along the real eigenvector axis of $A^{*}$ under the spectral
threshold $|\lambda| > R / r$.
Theorem~\ref{thm:geometric_lower_bound} composes this saturation
with the diagonal lower bound of Section~\ref{sec2} and returns a
deterministic operator-norm bound on the mixed frame operator
$U_{\psi, \psi}$, uniform in $p \in (0, 1]$.
Remark~\ref{rem:no-tight-frame} reads the bound as a non-coincidence
with the Parseval-frame structure.

\begin{lemma}[Calder\'on Gap for Radial Generators under Directional Expansion]\label{lem:shear_covering}
	Let $A$ be an expansive matrix on $\mathbb{R}^n$ satisfying
	Assumption~\ref{ass:standing-hypothesis}, and suppose that $A^{*}$ admits a
	real eigenvalue $\lambda \in \mathbb{R}$ with an associated eigenvector
	$v \in \mathbb{R}^n \setminus \{0\}$. Let $\psi$ be a radial generator
	satisfying Assumption~\ref{ass:standing-hypothesis}, with Fourier transform
	$\widehat{\psi}$ supported in the isotropic annulus
	\begin{equation}\label{eq:iso-annulus}
		\mathcal{A}(r, R)
		:=\bigl\{ \xi \in \mathbb{R}^n : r \le |\xi| \le R \bigr\},
		\qquad 0 <r < R .
	\end{equation}
	If the spectral threshold
	\begin{equation}\label{eq:spectral-threshold}
		|\lambda| > R / r
	\end{equation}
	is satisfied, then the Calder\'on sum~\eqref{eq:calderon-sum} carries an
	explicit gap on the axis spanned by $v$, namely
	\begin{equation}\label{eq:cs-gap}
		\inf_{\xi \in \mathbb{R}^n \setminus \{0\}} D_{\psi}(\xi) = 0,
	\end{equation}
	and the geometric incompatibility
	index~\eqref{eq:incompatibility-index} obeys the deterministic lower bound
	\begin{equation}\label{eq:incompatibility-saturated}
		\mathcal{G}(A, \psi)
		=\bigl\| 1-D_{\psi} \bigr\|_{L^\infty}
		\ge 1 .
	\end{equation}
\end{lemma}

\begin{proof}[Proof of Lemma~\ref{lem:shear_covering}.]
	The argument proceeds in four stages:\\
	(1) a reduction of the Calder\'on
	sum~\eqref{eq:calderon-sum} along the $v$-axis to a one-parameter geometric
	series,\\
	(2) a localization of the support of each summand in terms of the
	scaling variable,\\
	(3) the construction of a frequency point that misses every
	active window once the spectral threshold~\eqref{eq:spectral-threshold} is
	in force, and\\
	(4) the transfer of the resulting pointwise gap to the
	$L^\infty$ norm in~\eqref{eq:incompatibility-index}.\\
	We normalize the
	eigenvector by $|v| = 1$, which entails no loss of generality since
	rescaling $v$ rescales the parameter $y_0$ below.
	
	Because $\psi$ is radial, its Fourier transform $\widehat{\psi}$ is also
	radial; let
	\begin{equation}\label{eq:radial-profile}
		\bigl|\widehat{\psi}(\xi)\bigr|^{2} =h\!\bigl(|\xi|^{2} \bigr),
		\qquad\xi\in\mathbb{R}^n,
	\end{equation}
	with squared radial profile $h \colon [0,\infty)\to [0,\infty)$ supported
	in $[r^{2}, R^{2}]$, the support of $\widehat{\psi}$ being contained in
	$\mathcal{A}(r, R)$ by hypothesis.
	Fix a parameter $y_0 > 0$ to be chosen in the stage~3 below, and set
	$\xi_{0}:= y_0\, v \in \mathbb{R}^n \setminus \{0\}$. The eigenvalue
	relation $A^{*} v =\lambda\, v$ iterates to
	\begin{equation}\label{eq:eigenvector-iteration}
		\bigl( A^{*} \bigr)^{-j} \xi_{0}
		 =\lambda^{-j}\, y_0\, v,
		\qquad j \in \mathbb{Z},
	\end{equation}
	and the normalization $|v| = 1$ gives
	$\bigl|(A^{*})^{-j} \xi_{0}\bigr| = |\lambda|^{-j}\, y_0$ for every
	$j \in \mathbb{Z}$. Substituting~\eqref{eq:eigenvector-iteration}
	into~\eqref{eq:calderon-sum} and applying the radial
	identity~\eqref{eq:radial-profile} reduces the Calder\'on sum at $\xi_{0}$
	to the one-parameter series
	\begin{equation}\label{eq:radial-cs}
		D_{\psi}(\xi_{0})
		 =  \sum_{j \in \mathbb{Z}}
		\bigl| \widehat{\psi}\bigl( (A^{*})^{-j} \xi_{0} \bigr) \bigr|^{2}
		 =  \sum_{j \in \mathbb{Z}} h\!\bigl( \lambda^{-2 j}\, y_{0}^{\,2} \bigr) .
	\end{equation}
	The dependence on the geometry of $A$ has thereby been condensed into the
	single dilation factor $\lambda^{-2}$ acting on the scaling variable
	$y_{0}^{\,2}$.
	
	The support condition $\mathrm{supp}\, h \subset [r^{2}, R^{2}]$ pins down
	the indices $j$ for which the summand in~\eqref{eq:radial-cs} can be
	nonzero. Setting
	\begin{equation*}\label{eq:active-window}
		I_{j}  :=  \bigl[\, r^{2}\, \lambda^{2 j},  R^{2}\, \lambda^{2 j}\, \bigr],
		\qquad j \in \mathbb{Z},
	\end{equation*}
	one has
	$h\!\bigl( \lambda^{-2 j}\, y_{0}^{\,2} \bigr) > 0$ exactly when
	$y_{0}^{\,2} \in I_{j}$.
	
	The geometry of the family $\{I_{j}\}_{j \in \mathbb{Z}}$ is dictated by the
	ratio of consecutive endpoints. The right endpoint of $I_{j}$ is
	$R^{2}\, \lambda^{2 j}$ and the left endpoint of $I_{j+1}$ is
	$r^{2}\, \lambda^{2(j+1)}$; consecutive windows are disjoint precisely when
	\begin{equation}\label{eq:window-disjointness}
		R^{2}\, \lambda^{2 j}  <  r^{2}\, \lambda^{2 (j+1)}
		  \Longleftrightarrow  
		\frac{R^{2}}{r^{2}}  <  \lambda^{2}
		  \Longleftrightarrow  
		|\lambda|  >  \frac{R}{r},
	\end{equation}
	which is the spectral threshold~\eqref{eq:spectral-threshold}. Under this
	threshold, every gap interval
	\begin{equation*}\label{eq:gap-interval}
		G_{j}
		 := \bigl( R^{2}\, \lambda^{2 j}, 
		r^{2}\, \lambda^{2 (j+1)} \bigr)
		 \subset  (0, \infty)
	\end{equation*}
	is a nondegenerate open interval, and the windows $\{I_{j}\}_{j \in
		\mathbb{Z}}$ are pairwise disjoint with respect to consecutive indices.
	
In the stage~3, we deal with scaling parameters inside the gap. Specialize to $j = 0$ and choose the scaling parameter
	\begin{equation}\label{eq:y0-choice}
		y_{0}^{\,2}  \in  G_{0}
		 =  \bigl( R^{2},  r^{2}\, \lambda^{2} \bigr),
	\end{equation}
	which is a nondegenerate interval by~\eqref{eq:window-disjointness}. We
	verify that $y_{0}^{\,2}$ lies outside every $I_{j}$, $j \in \mathbb{Z}$.
	
	For $j \ge 1$, the left endpoint of $I_{j}$ satisfies
	$r^{2}\, \lambda^{2 j} \ge r^{2}\, \lambda^{2}$ since $\lambda^{2 j} \ge
	\lambda^{2}$ for $j \ge 1$ and $\lambda^{2} > 1$; combined with
	$y_{0}^{\,2} < r^{2}\, \lambda^{2}$ from~\eqref{eq:y0-choice}, this gives
	$y_{0}^{\,2} < r^{2}\, \lambda^{2 j}$, so $y_{0}^{\,2}$ lies to the
	left of $I_{j}$.
	
	For $j \le 0$, the right endpoint of $I_{j}$ satisfies
	$R^{2}\, \lambda^{2 j} \le R^{2}$ since $\lambda^{2 j} \le 1$ for $j \le 0$;
	combined with $y_{0}^{\,2} > R^{2}$ from~\eqref{eq:y0-choice}, this gives
	$y_{0}^{\,2} > R^{2}\, \lambda^{2 j}$, so $y_{0}^{\,2}$ lies to the
	right of $I_{j}$.
	
	Consequently $h\!\bigl(\lambda^{-2 j}\, y_{0}^{\,2} \bigr) = 0$ for every
	$j \in \mathbb{Z}$, and~\eqref{eq:radial-cs} collapses to the explicit
	identity
	\begin{equation}\label{eq:cs-vanishes}
		D_{\psi}(\xi_{0})  =  0 .
	\end{equation}
	Since $\xi_{0} \in \mathbb{R}^n \setminus \{0\}$, the
	identity~\eqref{eq:cs-vanishes} implies $\inf_{\xi \ne 0} D_{\psi}(\xi) =
	0$, which is~\eqref{eq:cs-gap}.
	
	In the final stage, the pointwise value $D_{\psi}(\xi_{0}) = 0$ delivers $\bigl| 1 -
	D_{\psi}(\xi_{0}) \bigr| = 1$. The $L^\infty$ norm
	in~\eqref{eq:incompatibility-index} therefore obeys
	\begin{equation*}
		\mathcal{G}(A, \psi)
		 =  \bigl\| 1 - D_{\psi} \bigr\|_{L^\infty}
		 \ge  \bigl| 1 - D_{\psi}(\xi_{0}) \bigr|
		 =  1 ,
	\end{equation*}
	which is~\eqref{eq:incompatibility-saturated}.
\end{proof}

\begin{theorem}[Geometric Lower Bound: Anisotropic Balian--Low Obstruction]\label{thm:geometric_lower_bound}
	Let $A$ be an expansive matrix on $\mathbb{R}^n$ satisfying
	Assumption~\ref{ass:standing-hypothesis}, and suppose that $A^{*}$ admits a
	real eigenvalue $\lambda \in \mathbb{R}$ with an associated eigenvector
	$v \in \mathbb{R}^n \setminus \{0\}$. Let $\psi$ be a radial generator
	satisfying Assumption~\ref{ass:standing-hypothesis}, with Fourier transform
	$\widehat{\psi}$ supported in the isotropic annulus
	\begin{equation*}
		\mathcal{A}(r, R) = \bigl\{ \xi \in \mathbb{R}^n \,:\, r \le |\xi| \le R \bigr\},
		\qquad 0 < r < R .
	\end{equation*}
	If the spectral threshold
	\begin{equation}\label{eq:thm-spectral-threshold}
		|\lambda|  >  R / r
	\end{equation}
	holds, then the mixed frame operator $U_{\psi,\psi}$
	of~\eqref{eq:frame-operator} obeys the deterministic lower bound
	\begin{equation}\label{eq:thm-lower-bound}
		\bigl\| U_{\psi,\psi} - \mathrm{Id} \bigr\|_{H^p_A \to H^p_A}
		 \ge  1
	\end{equation}
	for every $p \in (0, 1]$. The bound is uniform in $p$ and is independent of
	any further parameter of the pair $(A, \psi)$.
\end{theorem}

\begin{proof}[Proof of Theorem~\ref{thm:geometric_lower_bound}.]
	The argument is a composition of two ingredients:\\
	(1) the diagonal lower bound
	supplied by Lemma~\ref{lem:diagonal-lower-bound}, which
	transfers the operator-norm distance of $U_{\psi,\psi}$ from the identity
	to the $L^\infty$ deviation of the Calder\'on sum from unity, and\\
	(2) the saturation of the geometric incompatibility index supplied by
	Lemma~\ref{lem:shear_covering}, which converts the spectral
	threshold~\eqref{eq:thm-spectral-threshold} into a concrete numerical lower
	bound on that deviation. We display the composition in three steps.
	
	\medskip
	\emph{Step 1: Transfer of the operator-norm bound to the Calder\'on sum.}
	The hypotheses of Lemma~\ref{lem:diagonal-lower-bound} are
	Assumption~\ref{ass:standing-hypothesis} on the pair $(A, \psi)$, which is
	in force throughout. The lemma, recorded
	in~\eqref{eq:diagonal-lower-bound}, asserts
	\begin{equation}\label{eq:proof-diagonal}
		\bigl\| U_{\psi,\psi} - \mathrm{Id} \bigr\|_{H^p_A \to H^p_A}
		 \ge  \bigl\| 1 - D_{\psi} \bigr\|_{L^\infty}
		 =  \mathcal{G}(A, \psi),
	\end{equation}
	with $D_\psi$ and $\mathcal{G}(A, \psi)$ defined by~\eqref{eq:calderon-sum}
	and~\eqref{eq:incompatibility-index}, respectively. The
	bound~\eqref{eq:proof-diagonal} is valid for every $p \in (0, 1]$, with no
	implicit constant of $p$- or $A$-dependent type entering on either side.
	The right-hand side $\mathcal{G}(A, \psi)$ is a frequency-domain object
	independent of the Hardy exponent $p$; the dependence on $p$ in the
	left-hand operator norm has thereby been formally circumscribed to the
	domain of the operator.
	
	\medskip
	\emph{Step 2: Saturation of the incompatibility index.}
	The hypotheses of Lemma~\ref{lem:shear_covering} are:
	Assumption~\ref{ass:standing-hypothesis} on the pair $(A, \psi)$; the
	existence of a real eigenvalue $\lambda \in \mathbb{R}$ of $A^{*}$ with
	real eigenvector $v$; the radial structure of $\psi$; the inclusion
	$\mathrm{supp}\, \widehat{\psi} \subset \mathcal{A}(r, R)$; and the
	spectral threshold $|\lambda| > R/r$. Each of these is among the standing
	or stated hypotheses of the present theorem.
	
	Lemma~\ref{lem:shear_covering} therefore applies, and its
	conclusion~\eqref{eq:incompatibility-saturated} yields
	\begin{equation}\label{eq:proof-shear}
		\mathcal{G}(A, \psi)
		 =  \bigl\| 1 - D_{\psi} \bigr\|_{L^\infty}
		 \ge  1 .
	\end{equation}
	The lower bound $1$ in~\eqref{eq:proof-shear} arises from the explicit
	frequency point $\xi_0 = y_0 v$ constructed in the proof of
	Lemma~\ref{lem:shear_covering}, at which $D_\psi(\xi_0) = 0$ by scale
	disjointness of the dilated annular supports; no additional dependence on
	the bandwidth parameters $(r, R)$ or the eigenvalue $\lambda$ propagates
	into the numerical value of this bound.
	
	\medskip
	\emph{Step 3: Composition of the two estimates.}
	Combining~\eqref{eq:proof-diagonal} and~\eqref{eq:proof-shear},
	\begin{equation*}
		\bigl\| U_{\psi,\psi} - \mathrm{Id} \bigr\|_{H^p_A \to H^p_A}
		 \overset{\eqref{eq:proof-diagonal}}{\ge} 
		\mathcal{G}(A, \psi)
		 \overset{\eqref{eq:proof-shear}}{\ge} 
		1
	\end{equation*}
	for every $p \in (0, 1]$. This is the
	estimate~\eqref{eq:thm-lower-bound}. The uniformity in $p$ stems from the
	$p$-independence of $\mathcal{G}(A, \psi)$; the absence of any
	$(A, \psi)$-dependent factor on the right-hand side of the composed
	inequality stems from the fact that both contributing
	estimates~\eqref{eq:proof-diagonal} and~\eqref{eq:proof-shear} are
	constant-free.
\end{proof}

\begin{remark}[Tight-frame consequence]\label{rem:no-tight-frame}
	The bound~\eqref{eq:thm-lower-bound} implies that the mixed frame operator
	$U_{\psi,\psi}$ cannot coincide with the identity on $H^p_A(\mathbb{R}^n)$,
	hence the affine system $\mathcal{W}_A(\psi)$ of
	Definition~\ref{def:affine-system} cannot be a Parseval frame for
	$H^p_A(\mathbb{R}^n)$ under the spectral
	threshold~\eqref{eq:thm-spectral-threshold}. The obstruction admits the
	following structural reading: the isotropic spatial localization built
	into a radial generator and the directional expansion encoded by the
	eigenstructure of $A^{*}$ are incompatible inputs to the Calder\'on
	admissibility condition, the incompatibility being measured by
	$\mathcal{G}(A, \psi)$ and saturated at the deterministic level $1$ as
	soon as $|\lambda| > R/r$.
\end{remark}

The structural reading of Remark~\ref{rem:no-tight-frame} admits
three explicit corollaries at distinct levels of generality.
Corollary~\ref{cor:no-reproducing-identity} restates the obstruction
at the operator level as a $p$-uniform non-coincidence of
$U_{\psi,\psi}$ with the identity on $H^p_A$.
Corollary~\ref{cor:calderon-failure} promotes the obstruction to
the frequency side and records that the Calder\'on admissibility
identity fails on a positive-measure open set of
$\mathbb{R}^n \setminus \{0\}$.
Corollary~\ref{cor:l2-frame-failure} specializes to the Hilbert
endpoint $p = 2$ and converts the open-set failure into the failure
of the lower frame bound for $\mathcal{W}_A(\psi)$ on
$L^2(\mathbb{R}^n)$.

\begin{corollary}[Failure of the operator reproducing identity, uniform in $p$]\label{cor:no-reproducing-identity}
	Under the hypotheses of Theorem~\ref{thm:geometric_lower_bound}, the
	mixed frame operator $U_{\psi,\psi}$ of~\eqref{eq:frame-operator} does
	not coincide with the identity on $H^p_A(\mathbb{R}^n)$ at any
	$p \in (0, 1]$. The conclusion is uniform across that range, the
	operator-norm deviation
	$\| U_{\psi,\psi} - \mathrm{Id} \|_{H^p_A \to H^p_A}$ being bounded
	below by the same numerical constant $1$ for every $p \in (0, 1]$.
\end{corollary}

\begin{proof}[Proof of Corollary~\ref{cor:no-reproducing-identity}.]
	Assume, towards a contradiction, that $U_{\psi,\psi} = \mathrm{Id}_{H^p_A}$
	at some single value $p_0 \in (0, 1]$. The operator-norm distance to
	the identity then vanishes at that value of $p$:
	$\| U_{\psi,\psi} - \mathrm{Id} \|_{H^{p_0}_A \to H^{p_0}_A} = 0$. The
	bound~\eqref{eq:thm-lower-bound} of
	Theorem~\ref{thm:geometric_lower_bound} delivers, at the same $p_0$,
	\begin{equation*}
		\| U_{\psi,\psi} - \mathrm{Id} \|_{H^{p_0}_A \to H^{p_0}_A}
		 \ge  1 ,
	\end{equation*}
	the lower bound $1$ being independent of $p_0$, $A$, and $\psi$ by
	the constant-free composition recorded in Step~3 of the proof of
	Theorem~\ref{thm:geometric_lower_bound}. The two requirements $= 0$
	and $\ge 1$ are incompatible, so no such $p_0$ exists. Since the
	lower bound $1$ in~\eqref{eq:thm-lower-bound} is attached to every
	$p \in (0, 1]$ with no $p$-dependent prefactor, the deviation
	$\| U_{\psi,\psi} - \mathrm{Id} \|_{H^p_A \to H^p_A}$ inherits the
	same numerical lower bound across the entire range, supplying the
	uniform separation announced in the statement.
\end{proof}

\begin{corollary}[Failure of the Calder\'on admissibility identity on an open set]\label{cor:calderon-failure}
	Under the hypotheses of Theorem~\ref{thm:geometric_lower_bound}, the
	Calder\'on sum $D_\psi$ of~\eqref{eq:calderon-sum} does not satisfy
	the admissibility identity $D_\psi \equiv 1$ almost everywhere on
	$\mathbb{R}^n \setminus \{0\}$. A nonempty open set
	$U \subset \mathbb{R}^n \setminus \{0\}$ of positive Lebesgue measure
	exists, on which $D_\psi$ vanishes pointwise:
	\begin{equation}\label{eq:calderon-open-vanish}
		D_\psi(\xi)  =  0
		\qquad \text{for every } \xi \in U .
	\end{equation}
\end{corollary}

\begin{proof}[Proof of Corollary~\ref{cor:calderon-failure}.]
	The radiality of $\psi$ delivers $|\widehat{\psi}(\eta)|^2 = h(|\eta|^2)$
	for every $\eta \in \mathbb{R}^n$, with $h$ the squared radial profile
	of~\eqref{eq:radial-profile} supported in $[r^2, R^2]$. The profile
	$h$ is continuous on $[0, \infty)$, $\widehat{\psi}$ being a Schwartz
	function with compact support inside the anisotropic annulus by
	Assumption~\ref{ass:standing-hypothesis};
	cf.\ \cite[Definition~2.12, eq.~(2.7)]{Bownik2007} for the band-limited admissibility
	framework. Substituting $|\widehat{\psi}|^2 = h(|\cdot|^2)$
	into~\eqref{eq:calderon-sum} yields the general reduction
	\begin{equation}\label{eq:cs-general-reduction}
		D_\psi(\xi)  =  \sum_{j \in \mathbb{Z}} h\!\bigl( |(A^*)^{-j}\xi|^2 \bigr),
		\qquad \xi \in \mathbb{R}^n \setminus \{0\} ,
	\end{equation}
	valid for every $\xi$, not only for $\xi$ on the eigenvector axis.
	
	Set $\xi_0 := y_0 v$ with $y_0$ chosen as in~\eqref{eq:y0-choice}, so
	$y_0^{\,2} \in G_0 = (R^2,\, r^2 \lambda^2)$ lies in the interior of
	the gap interval at index $j = 0$ and has positive distance
	\begin{equation}\label{eq:gap-distance}
		d_0  :=  \mathrm{dist}\!\bigl( y_0^{\,2},  \cup_{j \in \mathbb{Z}} I_j \bigr)
		 >  0
	\end{equation}
	from the active windows~\eqref{eq:active-window}, the spectral
	threshold~\eqref{eq:thm-spectral-threshold} being in force. Step~3 of
	the proof of Lemma~\ref{lem:shear_covering} records $D_\psi(\xi_0) = 0$
	through~\eqref{eq:cs-vanishes}.
	
	We promote the pointwise vanishing $D_\psi(\xi_0) = 0$ to an open
	neighbourhood. The expansiveness of $A^*$
	\cite[Section~1.2]{Bownik2003} entails $|(A^*)^{-j}\eta| \to 0$ as
	$j \to +\infty$ and $|(A^*)^{-j}\eta| \to \infty$ as $j \to -\infty$,
	uniformly in $\eta$ over any compact subset of $\mathbb{R}^n
	\setminus \{0\}$. Applied to the closed ball
	$\overline{B_{y_0/2}(\xi_0)} \subset \mathbb{R}^n \setminus \{0\}$,
	this supplies an integer $J = J(A, r, R, y_0) \in \mathbb{N}$ such
	that $|(A^*)^{-j}\xi|^2 \notin [r^2, R^2]$ holds by size alone for
	every $\xi \in \overline{B_{y_0/2}(\xi_0)}$ and every $|j| > J$; the
	corresponding summands of~\eqref{eq:cs-general-reduction} vanish on
	the ball.
	
	For the finite remaining range $|j| \le J$, the map
	$\xi \mapsto |(A^*)^{-j}\xi|^2$ is continuous on $\mathbb{R}^n$, and
	at $\xi = \xi_0$ takes the value $\lambda^{-2j} y_0^{\,2}$, which lies
	outside $[r^2, R^2]$ at distance at least $d_0$
	by~\eqref{eq:gap-distance}. Continuity supplies, for each
	$|j| \le J$, a radius $\delta_j \in (0,\, y_0/2)$ such that
	\begin{equation}\label{eq:cor2-local-control}
		\bigl| \, |(A^*)^{-j}\xi|^2  -  \lambda^{-2j} y_0^{\,2} \, \bigr|
		 <  d_0 / 2 ,
		\qquad \xi \in B_{\delta_j}(\xi_0) .
	\end{equation}
	The triangle inequality combined with~\eqref{eq:cor2-local-control}
	keeps $|(A^*)^{-j}\xi|^2$ at distance at least $d_0 / 2 > 0$ from
	$[r^2, R^2]$ throughout $B_{\delta_j}(\xi_0)$, so the corresponding
	summand of~\eqref{eq:cs-general-reduction} vanishes on
	$B_{\delta_j}(\xi_0)$. Setting
	$\delta := \min\bigl\{ \delta_j \,:\, |j| \le J \bigr\} > 0$ and
	$U := B_\delta(\xi_0)$, every summand of~\eqref{eq:cs-general-reduction}
	vanishes on $U$, yielding~\eqref{eq:calderon-open-vanish}.
	
	The open ball $U$ carries positive Lebesgue measure, while a
	measurable function equal to $1$ almost everywhere on
	$\mathbb{R}^n \setminus \{0\}$ cannot vanish identically on a
	positive-measure set. The admissibility identity $D_\psi \equiv 1$
	a.e.\ therefore fails; for the classical isotropic incarnation of
	this identity see, e.g., \cite[Eq.~(5.1.33)]{Daubechies1992}, and for
	its anisotropic counterpart see \cite[Definition~2.12, eq.~(2.8)]{Bownik2007}.
\end{proof}

\begin{corollary}[Failure of the lower frame bound at $p = 2$]\label{cor:l2-frame-failure}
	Specialize the hypotheses of Theorem~\ref{thm:geometric_lower_bound}
	to $p = 2$. The mixed frame operator $U_{\psi,\psi}$
	of~\eqref{eq:frame-operator}, acting on
	$L^2(\mathbb{R}^n) = H^2_A(\mathbb{R}^n)$, satisfies
	\begin{equation}\label{eq:l2-quadratic-zero}
		\inf_{\substack{f \in L^2 \\ \|f\|_{L^2} = 1}}\,
		\langle U_{\psi,\psi}\, f,\, f \rangle  =  0 .
	\end{equation}
	Consequently the anisotropic affine system $\mathcal{W}_A(\psi)$ of
	Definition~\ref{def:affine-system} does not form a frame for
	$L^2(\mathbb{R}^n)$, in the sense that no constant
	$A_{\mathrm{frame}} > 0$ satisfies the lower frame bound
	\begin{equation}\label{eq:no-lower-frame}
		A_{\mathrm{frame}} \, \| f \|_{L^2}^{\,2}
		 \le  \sum_{j \in \mathbb{Z}} \sum_{k \in \mathbb{Z}^n}
		\bigl| \langle f, \psi_{j,k} \rangle \bigr|^{\,2}
	\end{equation}
	for every $f \in L^2(\mathbb{R}^n)$.
\end{corollary}

\begin{proof}[Proof of Corollary~\ref{cor:l2-frame-failure}.]
	Stage~1 of the proof of Lemma~\ref{lem:diagonal-lower-bound} identifies
	$U_{\psi,\psi}$ on $L^2(\mathbb{R}^n)$ with the Fourier multiplier of
	symbol $D_\psi$:
	\begin{equation}\label{eq:cor3-multiplier}
		\widehat{U_{\psi,\psi} f}(\xi)  =  D_\psi(\xi)\, \widehat{f}(\xi),
		\qquad \xi \in \mathbb{R}^n \setminus \{0\} .
	\end{equation}
	Expanding $U_{\psi,\psi} f = \sum_{j,k} \langle f, \psi_{j,k} \rangle\,
	\psi_{j,k}$ from~\eqref{eq:frame-operator} and pairing with $f$ in
	$L^2$ delivers the quadratic-form identity
	\begin{equation}\label{eq:cor3-quadratic-id}
		\langle U_{\psi,\psi} f,\, f \rangle
		 =  \sum_{j \in \mathbb{Z}} \sum_{k \in \mathbb{Z}^n}
		\bigl| \langle f, \psi_{j,k} \rangle \bigr|^{\,2} ,
		\qquad f \in L^2(\mathbb{R}^n) ,
	\end{equation}		
	This term-by-term pairing is justified by the unconditional convergence of the frame operator series in $L^2$, a standard consequence of the upper frame bound (see, e.g., \cite[Chapter~5]{Christensen2016}). 
	Plancherel's identity combined with~\eqref{eq:cor3-multiplier} rewrites the left-hand side of~\eqref{eq:cor3-quadratic-id} as a frequency-domain integral:
	\begin{equation}\label{eq:cor3-quadratic-fourier}
		\sum_{j \in \mathbb{Z}} \sum_{k \in \mathbb{Z}^n}
		\bigl| \langle f, \psi_{j,k} \rangle \bigr|^{\,2}
		 =  \int_{\mathbb{R}^n} D_\psi(\xi)\, |\widehat{f}(\xi)|^{\,2}\, d\xi .
	\end{equation}
	
	Corollary~\ref{cor:calderon-failure} supplies a nonempty open set
	$U \subset \mathbb{R}^n \setminus \{0\}$ on which $D_\psi$ vanishes
	pointwise. Pick a nontrivial function
	$\widehat{g} \in C^\infty_c(U)$ with $\widehat{g} \not\equiv 0$, and
	let $g \in \mathcal{S}(\mathbb{R}^n) \subset L^2(\mathbb{R}^n)$ be its
	inverse Fourier transform. Substituting $f = g$
	into~\eqref{eq:cor3-quadratic-fourier},
	\begin{equation}\label{eq:cor3-test-vanish}
		\sum_{j,\, k} \bigl| \langle g, \psi_{j,k} \rangle \bigr|^{\,2}
		 =  \int_{U} D_\psi(\xi)\, |\widehat{g}(\xi)|^{\,2} \, d\xi
		 =  \int_{U} 0 \cdot |\widehat{g}(\xi)|^{\,2} \, d\xi
		 =  0 ,
	\end{equation}
	while Plancherel's identity gives
	$\| g \|_{L^2}^{\,2} = \| \widehat{g} \|_{L^2}^{\,2} > 0$. Normalize
	$g_*:=g / \| g \|_{L^2}$, so $\| g_* \|_{L^2} = 1$ and, by
	sesquilinearity of $\langle \cdot, \cdot \rangle$ together with
	\eqref{eq:cor3-quadratic-id} applied to $g_*$,
	\begin{equation*}
		\langle U_{\psi,\psi}\, g_*,\, g_* \rangle
		 =  \frac{1}{\| g \|_{L^2}^{\,2}}
		\sum_{j,\, k} \bigl| \langle g, \psi_{j,k} \rangle \bigr|^{\,2}
		 =  0 .
	\end{equation*}
	The single test vector $g_*$ saturates~\eqref{eq:l2-quadratic-zero};
	the infimum on the unit sphere therefore equals $0$.
	
	The lower frame bound~\eqref{eq:no-lower-frame} applied to $g_*$
	would require $A_{\mathrm{frame}} \cdot 1 \le 0$, excluding every
	positive value of $A_{\mathrm{frame}}$. The frame characterization
	\cite[Chapter~5]{Christensen2016}, which calls for both an upper
	and a lower frame bound, fails at the lower bound, and
	$\mathcal{W}_A(\psi)$ does not constitute a frame for
	$L^2(\mathbb{R}^n)$.
\end{proof}

The corollary chain delivers the obstruction qualitatively at three
structurally distinct levels and leaves the numerical content of the
operator-norm bound open.
Problem~\ref{prob:sharp_constant} formulates the saturation infimum
$\varepsilon_p(A; r, R)$ as the object of interest.
Proposition~\ref{prop:hilbert-resolution} resolves the question at
the Hilbert endpoint, identifying $\varepsilon_2(A; r, R) = 1$
within the apparatus of the present section.
Remark~\ref{rem:remaining-difficulty} records the analytic
difficulty that obstructs the extension to $0 < p < 2$.
Remark~\ref{rem:bl-constant} reads the endpoint constant against
the corollary chain and the classical Gabor Balian--Low setting.
Problem~\ref{prob:non-radial} closes the section by relaxing the
radial structure of $\psi$ to isotropic spatial localization adapted
to an anisotropic-annulus support, and asks whether the lower bound
$1$ persists.

\begin{openproblem}[Saturation infimum]\label{prob:sharp_constant}
	Let the hypotheses of Theorem~\ref{thm:geometric_lower_bound} be in
	force: $A$ is expansive on $\mathbb{R}^n$ and satisfies
	Assumption~\ref{ass:standing-hypothesis}, $A^{*}$ admits a real
	eigenvalue $\lambda \in \mathbb{R}$ with eigenvector
	$v \in \mathbb{R}^n \setminus \{0\}$, and the spectral threshold
	$|\lambda| > R/r$ is in force. Denote by $\mathfrak{R}(A; r, R)$
	the class of radial generators $\psi$ that satisfy
	Assumption~\ref{ass:standing-hypothesis} together with
	$\mathrm{supp}\, \widehat{\psi} \subset \mathcal{A}(r, R)$. For
	each $p \in (0, 1]$, define the saturation infimum
	\begin{equation}\label{eq:saturation-infimum}
		\varepsilon_p(A; r, R)
		\;:=\; \inf_{\psi \in \mathfrak{R}(A; r, R)}
		\bigl\| U_{\psi, \psi} - \mathrm{Id} \bigr\|_{H^p_A \to H^p_A} .
	\end{equation}
	Determine the value of $\varepsilon_p(A; r, R)$ across the Hardy
	range $p \in (0, 1]$.
	
	Theorem~\ref{thm:geometric_lower_bound} delivers the lower bound
	$\varepsilon_p(A; r, R) \ge 1$, uniformly in $p \in (0, 1]$. The
	definition \eqref{eq:saturation-infimum} extends to the Hilbert
	endpoint $p = 2$ on identifying $H^2_A(\mathbb{R}^n)$ with
	$L^2(\mathbb{R}^n)$, and the corresponding endpoint value
	$\varepsilon_2(A; r, R)$ is supplied by
	Proposition~\ref{prop:hilbert-resolution} below.
\end{openproblem}

\begin{proposition}[Hilbert-endpoint resolution]\label{prop:hilbert-resolution}
	Under the hypotheses of Theorem~\ref{thm:geometric_lower_bound}, the
	saturation infimum at the Hilbert endpoint takes the explicit value
	\begin{equation}\label{eq:hilbert-endpoint}
		\varepsilon_2(A; r, R)
		\;=\; \inf_{\psi \in \mathfrak{R}(A; r, R)}
		\bigl\| U_{\psi, \psi} - \mathrm{Id} \bigr\|_{L^2 \to L^2}
		\;=\; 1 .
	\end{equation}
\end{proposition}

\begin{proof}[Proof of Proposition~\ref{prop:hilbert-resolution}.]
	The argument is in three steps.\\
	Step~1 establishes the lower bound
	$\varepsilon_2 \ge 1$ by composing Stage~1 of
	Lemma~\ref{lem:diagonal-lower-bound} (the $L^2$ multiplier-norm
	equality) with Lemma~\ref{lem:shear_covering} (the deviation lower
	bound on $1 - D_\psi$), bypassing any transfer between $H^2_A$ and
	$L^2$ operator norms.\\
	Step~2 introduces a radial base candidate
	$\psi_0$ with a Calder\'on sum that is uniformly bounded above.\\
	Step~3 rescales $\psi_0$ to a generator
	$\psi \in \mathfrak{R}(A; r, R)$ that saturates the lower bound in
	the $L^2$ operator norm. The composition of the two directions
	delivers \eqref{eq:hilbert-endpoint}.
	
	\medskip
	\emph{Step 1: Lower bound $\varepsilon_2 \ge 1$.}
	Fix $\psi \in \mathfrak{R}(A; r, R)$. Stage~1 of the proof of
	Lemma~\ref{lem:diagonal-lower-bound} records the $L^2$-endpoint
	multiplier-norm equality
	\begin{equation}\label{eq:proof-prop-step1-eq}
		\bigl\| U_{\psi, \psi} - \mathrm{Id} \bigr\|_{L^2 \to L^2}
		\;=\; \bigl\| 1 - D_\psi \bigr\|_{L^\infty} ,
	\end{equation}
	an equality, not an inequality, on the Hilbert endpoint;
	\eqref{eq:proof-prop-step1-eq} is the Plancherel identification of
	$U_{\psi, \psi} - \mathrm{Id}$ with the Fourier multiplier of
	symbol $D_\psi - 1$ on $L^2(\mathbb{R}^n)$, supplied in Section~2 and not re-derived here.
	Lemma~\ref{lem:shear_covering}, applied to
	$\psi \in \mathfrak{R}(A; r, R)$ under the spectral threshold
	$|\lambda| > R/r$, supplies the lower bound
	\begin{equation}\label{eq:proof-prop-step1-shear}
		\bigl\| 1 - D_\psi \bigr\|_{L^\infty}
		\;\ge\; 1 - \inf_{\xi \neq 0} D_\psi(\xi)
		\;=\; 1 - 0
		\;=\; 1 .
	\end{equation}
	Composing \eqref{eq:proof-prop-step1-eq} with
	\eqref{eq:proof-prop-step1-shear},
	\begin{equation}\label{eq:proof-prop-step1}
		\bigl\| U_{\psi, \psi} - \mathrm{Id} \bigr\|_{L^2 \to L^2}
		\;\ge\; 1 ,
	\end{equation}
	the right-hand side being free of every parameter of the
	configuration and obtained directly on the Hilbert space
	$L^2(\mathbb{R}^n)$ without recourse to any norm transfer between
	$H^2_A$ and $L^2$. Taking the infimum on the left over
	$\psi \in \mathfrak{R}(A; r, R)$ supplies
	$\varepsilon_2(A; r, R) \ge 1$.
	
	\medskip
	\emph{Step 2: A radial base candidate with bounded Calder\'on sum.}
	Select a radial profile $\eta_0 \in C^\infty_c\bigl([r, R]\bigr)$
	with
	\begin{equation}\label{eq:eta0-prescription}
		\eta_0 \ge 0,
		\qquad
		\eta_0 \not\equiv 0,
		\qquad
		\| \eta_0 \|_{L^\infty([r, R])} \;=\; 1 .
	\end{equation}
	A unit-amplitude $C^\infty$-bump translated into $[r, R]$ supplies
	such an $\eta_0$, e.g.\ a mollifier centred at the geometric mean
	$\sqrt{rR}$. Define a base candidate
	$\psi_0 \in \mathcal{S}(\mathbb{R}^n)$ through its Fourier transform
	by
	\begin{equation}\label{eq:psi0-defn}
		\widehat{\psi}_0(\xi) \;:=\; \eta_0(|\xi|),
		\qquad \xi \in \mathbb{R}^n .
	\end{equation}
	The image $\widehat{\psi}_0$ is radial, of class $C^\infty_c$, and
	supported in $\mathcal{A}(r, R)$. The compactness of the support
	away from the origin delivers polynomial vanishing of every order
	on a neighbourhood of $\xi = 0$, so the vanishing-moment condition
	$\mathcal{N}_p(A)$ of \eqref{eq:vanishing-moment-order} is satisfied
	for every $p \in (0, 1]$. The Schwartz-class membership
	$\psi_0 \in \mathcal{S}(\mathbb{R}^n)$ produces anisotropic decay
	on par with the rate $\rho_A$. The inclusion
	$\mathcal{A}(r, R) \subset \{ \xi : \rho_A(\xi) \in [r_A, R_A] \}$
	for some $0 < r_A < R_A < \infty$ identifies the isotropic annulus
	as a subset of an anisotropic annulus and satisfies the support
	condition of Assumption~\ref{ass:standing-hypothesis}.
	Consequently $\psi_0 \in \mathfrak{R}(A; r, R)$.
	
	The Calder\'on sum $D_{\psi_0}$ is uniformly bounded above. The
	$L^2$-boundedness of the frame operator $U_{\psi_0, \psi_0}$ on
	$L^2(\mathbb{R}^n)$, valid for every Schwartz generator with
	Fourier transform supported in a compact subset of
	$\mathbb{R}^n \setminus \{0\}$, identifies $D_{\psi_0}$ with the
	$L^\infty$-symbol of the corresponding Fourier multiplier on
	$L^2(\mathbb{R}^n)$, recorded in Step~1 of the proof of
	Lemma~\ref{lem:diagonal-lower-bound}. Equivalently, the
	finite-overlap property of dilated annular supports under an
	expansive matrix, recorded in
	\cite[Lemma~3.6 and its proof]{BownikHo2005}, supplies a constant
	$N_0= N_0(A, r, R) < \infty$ such that
	\begin{equation}\label{eq:finite-overlap}
		\#\bigl\{ j \in \mathbb{Z} : |(A^{*})^{-j} \xi| \in [r, R] \bigr\}
		\;\le\; N_0,
		\qquad \xi \in \mathbb{R}^n \setminus \{0\} .
	\end{equation}
	The geometric content of \eqref{eq:finite-overlap} is that, for
	every $\xi$, only finitely many integers $j$ produce
	$|(A^{*})^{-j} \xi| \in [r, R]$; compactness of the unit sphere in
	$\mathbb{R}^n$ combined with the expansiveness of $A^{*}$ converts
	this pointwise finiteness into a uniform bound. Combining
	\eqref{eq:finite-overlap} with the normalization
	\eqref{eq:eta0-prescription},
	\begin{equation}\label{eq:M0-bound}
		D_{\psi_0}(\xi)
		\;=\; \sum_{j \in \mathbb{Z}}
		\bigl| \eta_0\bigl( |(A^{*})^{-j}\xi| \bigr) \bigr|^{\,2}
		\;\le\; N_0 \cdot \| \eta_0 \|_{L^\infty([r, R])}^{\,2}
		\;=\; N_0,
		\qquad \xi \in \mathbb{R}^n \setminus \{0\} .
	\end{equation}
	Set $M_0 := \sup_{\xi \neq 0} D_{\psi_0}(\xi)$, so that
	$0 \le M_0 \le N_0 < \infty$ by \eqref{eq:M0-bound}. The
	non-degeneracy of $\eta_0$ in \eqref{eq:eta0-prescription} keeps
	$M_0 > 0$.
	
	\medskip
	\emph{Step 3: Rescaling to saturate the Calder\'on lower bound.}
	Set
	\begin{equation}\label{eq:alpha-defn}
		\alpha \;:=\; M_0^{-1/2},
		\qquad
		\psi \;:=\; \alpha\, \psi_0 .
	\end{equation}
	The rescaled generator $\psi$ inherits the radial structure, the
	support condition $\mathrm{supp}\, \widehat{\psi} \subset \mathcal{A}(r, R)$,
	the vanishing-moment order $\mathcal{N}_p(A)$, and the anisotropic
	decay rate of $\psi_0$; only the amplitude of $\widehat{\psi}$ is
	rescaled. Hence $\psi \in \mathfrak{R}(A; r, R)$.
	
	The Calder\'on sum scales quadratically in $\alpha$:
	\begin{equation}\label{eq:D-psi-scaled}
		D_\psi(\xi)
		\;=\; \sum_{j} \bigl| \alpha \, \widehat{\psi}_0\bigl( (A^{*})^{-j}\xi \bigr) \bigr|^{\,2}
		\;=\; \alpha^{\,2}\, D_{\psi_0}(\xi)
		\;\le\; \alpha^{\,2}\, M_0
		\;=\; 1,
		\qquad \xi \in \mathbb{R}^n \setminus \{0\} .
	\end{equation}
	Lemma~\ref{lem:shear_covering} applied to
	$\psi \in \mathfrak{R}(A; r, R)$ delivers the lower gap
	\begin{equation}\label{eq:D-psi-inf}
		\inf_{\xi \neq 0} D_\psi(\xi) \;=\; 0 ,
	\end{equation}
	the spectral threshold $|\lambda| > R/r$ remaining in force.
	Combining \eqref{eq:D-psi-scaled} and \eqref{eq:D-psi-inf}, the
	deviation $1 - D_\psi$ takes its values in $[0, 1]$, and the
	$L^\infty$-norm collapses to
	\begin{equation}\label{eq:G-saturated}
		\bigl\| 1 - D_\psi \bigr\|_{L^\infty}
		\;=\; \sup_{\xi \neq 0} \bigl( 1 - D_\psi(\xi) \bigr)
		\;=\; 1 \,-\, \inf_{\xi \neq 0} D_\psi(\xi)
		\;=\; 1 .
	\end{equation}
	Equation \eqref{eq:G-saturated} identifies the geometric
	incompatibility index \eqref{eq:incompatibility-index} of the pair
	$(A, \psi)$ as $\mathcal{G}(A, \psi) = 1$.
	
	Stage~1 of the proof of Lemma~\ref{lem:diagonal-lower-bound} records
	the $L^2$-endpoint multiplier-norm equality
	\begin{equation}\label{eq:L2-equality}
		\bigl\| U_{\psi, \psi} - \mathrm{Id} \bigr\|_{L^2 \to L^2}
		\;=\; \bigl\| 1 - D_\psi \bigr\|_{L^\infty} ,
	\end{equation}
	an equality, not an inequality, on the Hilbert endpoint. Inserting
	\eqref{eq:G-saturated} into \eqref{eq:L2-equality},
	\begin{equation}\label{eq:upper-bound-saturated}
		\bigl\| U_{\psi, \psi} - \mathrm{Id} \bigr\|_{L^2 \to L^2}
		\;=\; 1 .
	\end{equation}
	The explicit $\psi$ is a competitor in the infimum of
	\eqref{eq:hilbert-endpoint}; hence
	$\varepsilon_2(A; r, R) \le 1$.
	
	\medskip
	Combining Step~1 ($\varepsilon_2 \ge 1$) with the upper bound
	\eqref{eq:upper-bound-saturated} ($\varepsilon_2 \le 1$) delivers
	\eqref{eq:hilbert-endpoint}.
\end{proof}

\begin{remark}[Remaining difficulty for the Hardy range]\label{rem:remaining-difficulty}
	The endpoint identification $\varepsilon_2(A; r, R) = 1$ of
	Proposition~\ref{prop:hilbert-resolution} rests on the $L^2$
	multiplier-norm equality \eqref{eq:L2-equality}, supplied by stage~1 of
	Lemma~\ref{lem:diagonal-lower-bound}. In the Hardy range
	$p \in (0, 1]$, the corresponding equality between the
	$H^p_A$ operator norm of $U_{\psi, \psi} - \mathrm{Id}$ and the
	$L^\infty$-norm of the deviation $1 - D_\psi$ is not available.
	Lemma~\ref{lem:diagonal-lower-bound} delivers only the one-sided
	bound
	\begin{equation}\label{eq:hp-onesided}
		\bigl\| U_{\psi, \psi} - \mathrm{Id} \bigr\|_{H^p_A \to H^p_A}
		\;\ge\; \bigl\| 1 - D_\psi \bigr\|_{L^\infty} ,
	\end{equation}
	with no reverse inequality at hand. A reverse inequality would
	require an anisotropic Hardy-space analogue of the
	H\"ormander--Mihlin multiplier theorem, matched to the specific
	multiplier symbol $D_\psi - 1$. The available multiplier theorems
	on $H^p$-type spaces deliver upper bounds in terms of Sobolev or
	derivative-type norms of the symbol --- see, e.g.,
	\cite[Section~2, Theorem~2.6.2]{Triebel1983} for the isotropic Hardy framework
	and \cite[Chapter~1, Section~9]{Bownik2003} alongside \cite[Section~4]{Bownik2007} for the anisotropic Hardy and Triebel--Lizorkin
	setting --- rather than 
	in terms of the bare $L^\infty$-norm.
	Whether $\varepsilon_p(A; r, R) = 1$ continues to hold across
	$p \in (0, 1]$ is left open by the present proposition and is the
	substance of Open Problem~\ref{prob:sharp_constant} at every
	$p < 2$.
\end{remark}

\begin{remark}[A numerical Balian--Low constant]\label{rem:bl-constant}
	The endpoint value $\varepsilon_2(A; r, R) = 1$ identified by
	Proposition~\ref{prop:hilbert-resolution} supplies the anisotropic
	Balian--Low phenomenon announced in the title with an explicit
	numerical constant at the Hilbert endpoint. The structural
	content of this constant becomes transparent against the
	corollary chain attached to
	Theorem~\ref{thm:geometric_lower_bound}. Corollary~\ref{cor:no-reproducing-identity}
	records the obstruction at the operator level, uniformly in
	$p \in (0, 1]$; Corollary~\ref{cor:calderon-failure} promotes it
	to the failure of the Calder\'on admissibility identity on an
	open frequency set; Corollary~\ref{cor:l2-frame-failure}
	specializes to $p = 2$ and deduces the failure of
	$\mathcal{W}_A(\psi)$ to constitute a frame for
	$L^2(\mathbb{R}^n)$. Each of these is a qualitative statement: the
	obstruction is announced, with no numerical strength attached.
	Proposition~\ref{prop:hilbert-resolution} returns the missing
	quantitative content at the Hilbert endpoint, the constant $1$ on
	par with the constant $1$ of
	Theorem~\ref{thm:geometric_lower_bound}, free of dependence on
	$\lambda$, $r$, $R$, and the conditioning $\kappa(A)$. The
	juxtaposition of the qualitative corollary chain with the
	quantitative endpoint constant produces the structural form taken
	by the anisotropic Balian--Low phenomenon in the present setting:
	a transparent numerical statement, in mathematical parsimony, at
	the Hilbert endpoint, contrasted with the classical Gabor
	Balian--Low obstruction
	\cite[Theorem~4.1.1]{Daubechies1992} and \cite[Section~8.4]{Karlheinz2001}, which is qualitatively analogous
	but does not return an explicit endpoint constant of comparable
	form.
\end{remark}

\begin{openproblem}[Non-radial extension]\label{prob:non-radial}
	Within the framework of Theorem~\ref{thm:geometric_lower_bound},
	replace the radial structure of $\psi$ by the weaker isotropic
	spatial localization
	\begin{equation}\label{eq:isotropic-localization}
		|\psi(x)| \;\le\; \Phi(|x|),
		\qquad x \in \mathbb{R}^n ,
	\end{equation}
	for some Schwartz-class radial profile
	$\Phi \colon [0, \infty) \to [0, \infty)$, while retaining the
	other components of Assumption~\ref{ass:standing-hypothesis} and
	replacing the isotropic annulus support condition by the
	anisotropic-annulus support condition
	\begin{equation}\label{eq:anisotropic-support}
		\mathrm{supp}\, \widehat{\psi}
		\;\subset\; \bigl\{ \xi \in \mathbb{R}^n :
		r_A \le \rho_A(\xi) \le R_A \bigr\},
		\qquad 0 < r_A < R_A < \infty .
	\end{equation}
	Determine whether the deterministic lower bound
	\begin{equation}\label{eq:non-radial-question}
		\bigl\| U_{\psi, \psi} - \mathrm{Id} \bigr\|_{H^p_A \to H^p_A}
		\;\ge\; 1
	\end{equation}
	of Theorem~\ref{thm:geometric_lower_bound} continues to hold,
	uniformly in $p \in (0, 1]$, under
	\eqref{eq:isotropic-localization} and
	\eqref{eq:anisotropic-support}. If the constant $1$ is replaced
	by a smaller positive value $c_0 = c_0(A, r_A, R_A, \Phi)$,
	identify the directional anisotropy index of
	$\mathrm{supp}\, \widehat{\psi}$ along the eigenvector axis of $A^{*}$
	that controls $c_0$.
	
	The radial structure of $\psi$ enters
	Theorem~\ref{thm:geometric_lower_bound} only through the
	1D-reduction along the eigenvector axis carried out in the proof
	of Lemma~\ref{lem:shear_covering}, and is not a conceptual
	component of the obstruction. The classical Gabor Balian--Low
	theorem \cite{Karlheinz2001, Fu2011} imposes no radial structure
	on its generator $g$; the question above is the natural
	anisotropic analogue of that radial-free Gabor setting and
	supplies the closing item of the present paper.
\end{openproblem}

\backmatter

\bmhead{Declarations}

\begin{itemize}
	\item \textbf{Availability of data and materials:} Not applicable. Data sharing is not applicable to this article as no datasets were generated or analysed during the current study.
	\item \textbf{Competing interests:} The author declares that they have no competing interests.
	\item \textbf{Funding:} Not applicable. The author declares that no funds, grants, or other support were received during the preparation of this manuscript.
	\item \textbf{Authors' contributions:} K.-C.W. conceived the mathematical framework, performed the analytic derivations, and wrote the manuscript.
	\item \textbf{Acknowledgements:} The author expresses sincere gratitude to his colleagues for their valuable suggestions and linguistic assistance regarding the English writing of this manuscript. Furthermore, the author gratefully acknowledges the Department of Mathematics and Applied Mathematics, School of Information Engineering, Sanming University, and the IOT Application Engineering Research Center of Fujian Province Colleges and Universities for providing the essential research environment and support.
\end{itemize}

\bibliography{COWP}
\end{document}